\documentclass[9pt,amsfonts]{article}
\usepackage{graphicx,amssymb,eucal,mathrsfs}
\usepackage{latexsym,amsmath,amscd,epsfig,epic,eepic,xcolor}
\usepackage[all]{xy}
\usepackage{wasysym,ntheorem}
\usepackage{epsf,graphicx,pgf,etex,amscd,pstricks}

\usepackage[hypertexnames=false,backref=page,pdftex,
 	pdfpagemode=UseNone,
 	breaklinks=true,
 	extension=pdf,
 	colorlinks=true,
 	linkcolor=blue,
 	citecolor=red,
 	urlcolor=blue,
 ]{hyperref}

\newcommand{\Z}{\mathbb{Z}}
\newcommand{\R}{\mathbb{R}}

\newcommand{\C}{\mathbb{C}}

\newcommand{\PUD}{{\rm PU}(2,1)}
\newcommand{\HC}{\mathbb{H}_{\mathbb{C}}^{2}}
\newcommand{\ssg}{\overline{\Gamma_{1}}}
\newcommand{\sgr}{\Gamma_{1}}

\newtheorem{theorem}{{\bf Theorem}}

\newtheorem{prop}[theorem]{Proposition}

\newtheorem{rem}{Remark}

\newtheorem{defi}{Definition}
\newtheorem{exam}{Example}

\theoremnumbering{Alph}
\newtheorem*{main}{{\bf  Main Theorem}}

\title{Some noncoherent, nonpositively curved K\"ahler groups}
\author{Pierre Py}
\date{October 2014}

\begin{document}

\maketitle

\begin{abstract} If $\Gamma$ is any nonuniform lattice in the group ${\rm PU}(2,1)$, let $\overline{\Gamma}$ be the quotient of $\Gamma$ obtained by filling the cusps of $\Gamma$ (i.e.\ killing the center of parabolic subgroups). Assuming that such a lattice $\Gamma$ has positive first Betti number, we prove that for any sufficiently deep subgroup of finite index $\Gamma_{1} < \Gamma$, the group $\overline{\Gamma_{1}}$ is noncoherent. The proof relies on previous work of M.~Kapovich as well as of C.~Hummel and V.~Schroeder.  
\end{abstract}

\tableofcontents

\newpage


\section{Introduction}

A group $G$ is called {\it coherent} if every finitely generated subgroup of $G$ is finitely presented. This notion has been studied within various classes of groups and has a long history. For instance, it is easy to see that the fundamental group of a closed orientable surface is coherent, and Scott proved that the fundamental group of any $3$-manifold is also coherent~\cite{scott}. On the other hand, it has been known for a long time that the direct product of two non-Abelian free groups is not coherent, see~\cite{grunewald}. This implies for instance that ${\rm SL}_{n}(\Z)$ is not coherent if $n\ge 4$ whereas the coherence of ${\rm SL}_{3}(\Z)$ is an old open problem, see~\cite{wall}. For other examples of incoherent groups the reader can consult Wise's article~\cite{wise2011}.

In~\cite{kgt13}, Kapovich conjectured that any lattice $\Gamma$ in a semisimple Lie group not locally isomorphic to ${\rm SL}_{2}(\R)$ or ${\rm SL}_{2}(\C)$ is not coherent (the conjecture for lattices in the isometry group of real hyperbolic spaces of dimension at least $4$ is due to Wise). Kapovich proved this conjecture for many rank $1$ lattices, see~\cite{kgt13} for a precise statement as well as~\cite{bm,kp,kpv,p1994} for earlier results which motivated the conjecture. He proved in particular the following statement: 
\begin{center}
{\it If $\Gamma$ is a cocompact lattice in the isometry group of the complex hyperbolic plane and if ${\rm vb}_{1}(\Gamma)$ is positive, than $\Gamma$ is noncoherent}. 
\end{center}
Here ${\rm vb}_{1}(\Gamma)$ denotes the virtual first Betti number of $\Gamma$, i.e.\ the supremum of the first Betti numbers of finite index subgroups of $\Gamma$. This result implies in particular that any cocompact arithmetic lattice of the simplest type in the group ${\rm PU}(n,1)$ of holomorphic isometries of the complex hyperbolic space of dimension $n$ is noncoherent.  

In this note we will see that the arguments of Kapovich can be used to give other examples of aspherical complex surfaces with noncoherent fundamental groups. The fundamental groups of the surfaces we will be dealing with are closely related to, but different from, complex hyberpolic lattices (and hence we will say nothing new about Kapovich's conjecture). Before stating our result, we recall classical facts concerning non-uniform lattices in the group $\PUD$ and their parabolic subgroups. We will use freely the notions of elliptic, parabolic and hyperbolic isometries for elements of $\PUD$. For a definition of these notions in the context of ${\rm CAT}(0)$ spaces, which applies in particular to symmetric spaces of noncompact type, we refer the reader to~\cite[II.6]{bh}.

So let $\Gamma\subset \PUD$ be a non-uniform lattice. Let $\xi$ be a point in the boundary of the complex hyperbolic plane $\HC$ and let $H_{\xi}$ be a fixed horopshere centered at $\xi$ inside $\HC$. Recall that $H_{\xi}$ can be identified with the $3$-dimensional Heisenberg group, that we will simply denote by $N$. Under this identification, and using the embedding
$$H_{\xi}\simeq N \hookrightarrow \HC,$$
\noindent
the metric of $\HC$ induces a left-invariant metric on $N$. The isometry group ${\rm Iso}(H_{\xi})$ of this metric is isomorphic to a semidirect product of the form $N \rtimes K$ where $K$ is a compact group. The point $\xi$ is called a parabolic point for $\Gamma$ if the stabilizer $\Gamma (\xi)$ of $\xi$ inside $\Gamma$ contains a parabolic isometry. In this case, any element in $\Gamma (\xi)-\{id\}$ is parabolic or elliptic, the group $\Gamma (\xi)$ is a lattice in the group ${\rm Iso}(H_{\xi})$ and is called a parabolic subgroup of $\Gamma$. For all of this, we refer the reader to~\cite{gr,hs}. It is known that the intersection of $\Gamma (\xi)$ with the normal subgroup $N\triangleleft {\rm Iso}(H_{\xi})$ has finite index in $\Gamma (\xi)$~\cite{auslander}. 

\begin{defi}\label{nicecusps}
We will say that the lattice $\Gamma$ has {\it nice cusps} if each parabolic subgroup $\Gamma ( \xi)$ is actually contained in the normal subgroup $N\triangleleft {\rm Iso}(H_{\xi})$ (following the terminology of~\cite{hummel1998}, we also say that the parabolic isometries of $\Gamma$ have no rotational part). 
\end{defi}

It is well-known that any nonuniform lattice $\Gamma < \PUD$ has a finite index subgroup with nice cusps; the reader will find a proof of this fact in~\cite{hummel1998} for instance. If such a lattice $\Gamma$ is torsion-free and has nice cusps, the ends of the manifold $\HC/\Gamma$ are diffeomorphic to the product of an interval by a nilmanifold. In this case we will denote by $\overline{\Gamma}$ the quotient of $\Gamma$ by the normal subgroup generated by the centers of all parabolic subgroups of $\Gamma$. The group $\overline{\Gamma}$ is the fundamental group of the natural {\it toroidal compactification} of the open complex surface $\HC/\Gamma$. We will spend some time describing this compactification in section~\ref{rappelcusp}. We will sometime informally refer to the group $\overline{\Gamma}$ as the {\it filling} of the lattice $\Gamma$. We can now state the:

\begin{main}\label{nc} Let $\Gamma < \PUD$ be a torsion-free nonuniform lattice with nice cusps. Assume that $b_{1}(\Gamma)$ is positive. Then, there exists a finite set $\mathscr{F}\subset \Gamma$ of parabolic isometries such that for any finite index normal subgroup $\Gamma_{1}\triangleleft \Gamma$ with $\Gamma_{1}\cap \mathscr{F}=\emptyset$, the group $\overline{\Gamma_{1}}$ is not coherent.  
\end{main}

As the reader will see, the strategy of the proof is very similar to the one used by Kapovich to prove the result mentioned earlier about cocompact lattices in $\PUD$. The main new ingredients we will need are a result of Hummel and Schroeder saying that the filling $\overline{\Gamma_{1}}$ appearing in the theorem is the fundamental group of a nonpositively curved Riemannian manifold, and Poincar\'e's reducibility theorem for Abelian varieties. 

Let us also mention that there are several examples of nonuniform lattices in $\PUD$ having positive virtual first Betti number, and hence having finite index subgroups satisfying the hypothesis of the theorem. This is the case for instance for the group ${\rm PU}(2,1,\mathbb{Z}[i])$~\cite{wallach}. More examples are provided by the lattices constructed by Deligne and Mostow~\cite{dm}; indeed among their lattices one can find nonunifom lattices in $\PUD$ for which the corresponding complex hyperbolic orbifold admits a holomorphic map to a hyperbolic Riemann surface and thus has positive first Betti number. We refer the reader to~\cite{deraux} for a study of certain holomorphic maps between Deligne-Mostow quotients.   

The text is organized as follows. In section~\ref{rappelcusp}, we recall classical facts concerning the compactification of finite volume quotients of the complex hyperbolic space and in particular we recall the results from~\cite{hs}. In section~\ref{rappelstrategie}, we describe Kapovich's strategy to study the coherence of fundamental groups of aspherical K\"ahler surfaces with positive first Betti number. Finally, we conclude the proof in section~\ref{proof}, using the description of Abelian subgroups of rank $2$ of the fillings of nonuniform lattices in $\PUD$, as well as Poincar\'e's reducibility theorem.

I would like to thank Yves de Cornulier for his comments on a preliminary version of this text. 


\section{Cusp closing after Hummel and Schroeder}\label{rappelcusp}

 We first state some classical results concerning cusps of complex hyperbolic manifolds of finite volume and their compactification by Abelian varieties. All the results that we will discuss in this section remain true for non-uniform lattices in ${\rm PU}(n,1)$ for any $n\ge 2$; however we will only state them in the case $n=2$ since this is the case we are dealing with in this article. These results are well-known and can be found for instance in~\cite{hs} or~\cite{mok}. One important point here is that this discussion is independent of the arithmeticity of the lattices under consideration. We follow the presentation and notations from~\cite{hs}.

We recall briefly the definition of the complex hyperbolic plane $\HC$; see~\cite{gold} for more details. We consider the vector space $\C^{3}$ endowed with the Hermitian form defined by:
$$\langle z,w\rangle=z_{1}\overline{w_{1}}+z_{2}\overline{w_{2}}-z_{3}\overline{w_{3}}.$$
The space $\HC$ is the open subset of $\mathbb{P}(\C^{3})$ made of lines which are negative for the form $\langle \cdot , \cdot \rangle$; its boundary $\partial \HC$ is the subet of $\mathbb{P}(\C^{3})$ made of isotropic lines. The space $\HC$ carries a $\PUD$-invariant K\"ahler metric of negative curvature for which the visual boundary is naturally identified with $\partial \HC$. We assume that the metric is normalized to have constant holomorphic sectional curvature equal to $-4$.

Let $\xi \in \partial \HC$ be a point in the boundary of the complex hyperbolic plane $\HC$. Let $\phi^{t}$ be the gradient flow of the negative of any Busemann function associated to $\xi$. Concretely, for $p\in \HC$ and $t\in \R$ the point $\phi^{t}(p)$ is the point at distance $t$ from $p$ on the oriented geodesic going from $p$ to $\xi$. Finally, let $N$ be the unipotent radical of the stabilizer of $\xi$ in the group ${\rm PU}(2,1)$; this is a simply connected nilpotent Lie group isomorphic to the Heisenberg group. A down-to-earth description of $N$ can be obtained as follows. Pick a basis $(v_{1}, v_{2}, v_{3})$ of $\C^{3}$ such that:
\begin{enumerate}
\item the vector $v_{1}$ spans the line $\xi$,
\item  the vectors $v_{1}$ and $v_{2}$ are isotropic and satisfy $\langle v_{1}, v_{2}\rangle =1$,
\item the vector $v_{3}$ satisfy $\langle v_{3},v_{i}\rangle=\delta_{i3}$ for $i=1, 2, 3$.  
\end{enumerate}

Representing linear transformations of $\C^{3}$ by their matrices in the basis $(v_{1}, v_{2}, v_{3})$, one checks that the group $N$ can be identified with the group of $3 \times 3$ matrices with complex entries of the form:

\begin{equation}\label{model}
\left( \begin{array}{ccc}
1 & \beta & -\overline{v} \\
0 & 1 & 0 \\
0 & v & 1 \\
\end{array}\right)\end{equation}
where $\beta$ and $v$ are complex numbers satisfying $2 Re (\beta)+\vert v \vert^{2}=0$. Hence for $s\in \R$ and $v\in \C$ we will denote by $g(s,v)$ the matrix above with $\beta=\frac{-\vert v \vert^{2}}{2}+is$, i.e.
\begin{equation}\label{mdc}
g(v,s)=\left( \begin{array}{ccc}
1 & \frac{-\vert v \vert^{2}}{2}+is & -\overline{v} \\
0 & 1 & 0 \\
0 & v & 1 \\
\end{array}\right).
\end{equation}
The center of $N$ is the group of matrices of the form $g(0,s)$ for $s\in \R$. We will denote by $\frak{n}$ the Lie algebra of $N$. 

We now fix a base point $o$ in $\HC$. The map 
\begin{equation}\label{iid}
\R \times N \to \HC
\end{equation}
sending $(t,n)$ to $\phi^{t}(n\cdot o)$ is a diffeomorphism. Let $\mu$ be the scalar product on $\frak{n}$ induced by the embedding $n\mapsto n\cdot o$, let $\frak{z}$ the Lie algebra of the center of $N$ and $\frak{r}$ the orthogonal of $\frak{z}$ for $\mu$. If $a$ and $b$ are positive real numbers, let $\mu_{a,b}$ be the scalar product on $\frak{n}$ which also makes $\frak{z}$ and $\frak{r}$ orthogonal, coincides with $a^{2}\mu$ on $\frak{r} \times \frak{r}$ and with $b^{2}\mu$ on $\frak{z}\times \frak{z}$. The pull-back of the metric of $\HC$ to $\R \times N$ has the form:
$$dt^{2}\oplus \mu_{e^{-t},e^{-2t}}.$$
In this expression, the scalar product $\mu_{e^{-t},e^{-2t}}$ on $\frak{n}$ is identified with a left-invariant metric on $N$. 

\begin{exam}\label{pointb} If one uses again the coordinates associated to a basis $(v_{1}, v_{2}, v_{3})$ as above and if $o$ is the point of $\HC$ whose homogeneous coordinates are $[-\frac{1}{2}: 1: 0]$, then $\frak{r}$ is the subspace of $\frak{n}$ of matrices of the form 
\begin{equation}
\left(\begin{array}{ccc}
0 & 0 & -\overline{v} \\
0 & 0 & 0 \\
0 & v & 0 \\
\end{array}\right), 
\end{equation}
and $\mu$ is the metric $ds^{2}+\vert v \vert^{2}$. 
\end{exam}

The following proposition describes the structure of the quotient of $\HC$ (or of a horoball of $\HC$ centered at $\xi$) by a cyclic subgroup of the center of $N$. Although this proposition is well-known and appears for instance in~\cite{mok}, we will explain its proof for the reader's convenience. 

\begin{prop}\label{modelcusp} Let $\ell$ be a positive real number. The quotient $\R \times \langle g(0,\ell)\rangle \backslash N$ is biholomorphic to the open subset $O$ of $\C^{2}$ defined as follows:
$$O=\{(v,w)\in \C^{2}, 0 < \vert w\vert < e^{\frac{-\pi \vert v \vert^{2}}{\ell}}\}$$
Under this identification, the projections of horoballs centered at $\xi$ in $\R \times \langle g(0,\ell)\rangle \backslash N$ corresponds to subsets of the form $O_{C}:=\{(v,w)\in \C^{2}, 0<\vert w \vert < Ce^{\frac{-\pi \vert v \vert^{2}}{\ell}}\}$ for some constant $C\in (0,1)$.  
\end{prop}

As the reader will see, the foliation defined by $v={\rm cst}$ in the coordinate system $(v,w)$ above has an intrinsic meaning: if one lifts this foliation to $\HC$, then the leaf through a point $p$ is the complex geodesic containing the real geodesic $[p,\xi]$. 

{\it Proof.} Let $Z\in \frak{z}$ be a vector of norm $1$ for $\mu$, identified with a left-invariant vector field on $N$ and hence with a vector field on $\R \times N$, and let $T$ be the vector field $\frac{\partial}{\partial t}$. Up to replacing $Z$ by $-Z$, the almost complex structure $J$ on $\R \times N$ coming from the identification~\eqref{iid} satisfies: 
\begin{enumerate}
\item $J(Z)=e^{-2t}T$,
\item $J(T)=-e^{2t}Z$,
\item the subbundle of the tangent bundle of $\R \times N$ generated by left-invariant vector fields in $\frak{r}$ is $J$-invariant. 
\end{enumerate}
The fact that the left-invariant subbundle of the tangent bundle associated to $\frak{r}$ is $J$-invariant implies that the Abelian group $N/[N,N]\simeq \R^{2}$ naturally carries an almost complex structure $\underline{J}$ which is invariant by translation. This means that $N/[N,N]$ can be identified with $\C$ in such a way that the projection 
$$\begin{array}{ccc}
\R \times N & \to & N/[N,N] \\
(t,n) & \mapsto & \pi(n)\\
\end{array}$$
is holomorphic, where $\pi : N \to N/[N, N]$ is the natural projection. Using our down-to-earth model~\eqref{model}, the map $\pi$ is simply the map taking $g(v,s)$ to $v$. 

To simplify the calculations, we will now assume that the base point $o$ is the point $[-\frac{1}{2}:1:0]$ considered in Example~\ref{pointb}. Let $ L : \R \times N \to \C$ be the function defined by $$L(t,g(v,s))=-e^{2t}-\vert v \vert^{2}+2is.$$ Using our previous description of the complex structure on $\R \times N$, and the description of the subspace $\frak{r}$ given in Example~\ref{pointb}, one checks that the map $L$ is holomorphic. When $v\in \C$ is fixed, the image of $L(\cdot , g(v,\cdot ))$ equals the open set $\{z\in \C, Re(z)<-\vert v \vert^{2}\}$, hence for fixed $v$ and varying $s$ and $t$, the complex number $$e^{\frac{\pi L(t,g(v,s))}{\ell}}$$ ranges in the punctured disc
$$\{ z\in \C, 0< \vert z \vert  < e^{\frac{-\pi \vert v \vert^{2}}{\ell}}\}.$$
Thus, the image of the map $$\varphi:=(\pi, e^{\frac{\pi L}{\ell}}) : \R \times N \to \C^{2}$$ 
is exactly the open set $O$ appearing in the statement of the proposition. This map descends to $\R \times \langle g(0,\ell)\rangle \backslash N$ and gives a holomorphic diffeomorphism between $\R \times \langle g(0,\ell)\rangle \backslash N$ and $O$.

Consider now the Busemann function on $\HC \simeq \R \times N $ given by $B(t,g(v,s))=-t$. One has $e^{\frac{-\pi}{\ell}(\vert v \vert^{2} +e^{-2B})}=\vert w \vert$. Hence subsets of the form $\{B\le {\rm cst}\}$ project via $\varphi$ onto subsets of the form $\vert w \vert \cdot e^{\frac{\pi \vert v \vert^{2}}{\ell}}\le {\rm cst}$. This proves the second point of the proposition.\hfill $\Box$

Consider now a lattice $\Lambda < N$. Let $\ell>0$ be such that the center $Z(\Lambda)$ of $\Lambda$ is generated by $g(0,\ell)$. Identify the quotient $\HC/Z(\Lambda)$ with the open set $O\subset \C^{2}$ appearing in Proposition~\ref{modelcusp}. The group $\Lambda/Z(\Lambda)$ acts on $O$, the projection of an element $g(v_{0},s_{0})$ acting by:
$$(v,w)\mapsto (v+v_{0},e^{\frac{2i\pi s_{0}-\pi \vert v_{0}\vert^{2}-2\pi v\cdot \overline{v_{0}}}{\ell}}w).$$  
The group $\Lambda/Z(\Lambda)$ projects injectively into $N/[N,N]\simeq \C$, its image is a lattice denoted by $A$. The quotient $O_{\Lambda}$ of $O$ by this action admits a submersion onto the elliptic curve $\C/A$, whose fibers are punctured discs. One can complete $O_{\Lambda}$ by considering the open set $\widetilde{O}=\{(v,w), \vert w \vert < e^{\frac{-\pi \vert v\vert^{2}}{\ell}}\}$. The action of $\Lambda/Z(\Lambda)$ on $O$ extends to an action on $\widetilde{O}$; the quotient $\widetilde{O}_{\Lambda}$ of $\widetilde{O}$ by this action is a disc bundle over the elliptic curve $\C/A$. One can do the same construction replacing $\HC$ at the beginning by a horoball and replacing the open set $O$ by the open set $O_{C}$ appearing in Proposition~\ref{modelcusp} for a suitable constant $C$.

Untill the end of Section~\ref{rappelcusp}, all the lattices of $\PUD$ that we consider are assumed to be torsion-free. Let now $\Gamma < \PUD$ be a nonuniform lattice with nice cusps. One knows that $\HC/\Gamma$ is diffeomorphic to the union of a compact submanifold with boundary, and finitely many cusps, which are, by definition, subsets of the form $[t_{i}, +\infty)\times \Lambda_{i}\backslash N$ for some parabolic subgroups $\Lambda_{i}$. By compactifying each cusp by the process described above, one obtains a compact complex surface that we will denote by $X_{\Gamma}$, which is the disjoint union of $\HC/\Gamma$ and finitely many elliptic curves. This compactification is canonical. Moreover, the fundamental group of $X_{\Gamma}$ is naturally isomorphic to the group $\overline{\Gamma}$ defined in the introduction. 

We now state some of the main results from~\cite{hs}:

\begin{enumerate}
\item First the complex surface $X_{\Gamma}$ is K\"ahler, see Theorem 7 in~\cite{hs}. 
\item Second, if $\Gamma$ is fixed, there exists a finite set $\mathscr{F}$ of parabolic isometries of $\Gamma$ such that if $\Gamma_{1}$ is a finite index normal subgroup of $\Gamma$ whose intersection with $\mathscr{F}$ is empty, then $X_{\Gamma_{1}}$ admits a Riemannian metric $h$ of nonpositive curvature, which has moreover the following property. The sectional curvature of $h$ along any $2$-plane $P\subset T_{x}X_{\Gamma_{1}}$ is negative is $x$ does not lie on one of the compactifying elliptic curves; these elliptic curves are flat and totally geodesic. See~\cite[\S 3]{hs}, in particular Proposition 3.3 and Remark 1 on pages 293-294.  
\end{enumerate}

\begin{rem} As for the original surface $X_{\Gamma}$, it always carries the structure of a nonpositively curved {\it orbifold}, see~\cite{hs}. 
\end{rem}

\begin{rem} If $\Gamma_{1}$ is a finite index normal subgroup of $\Gamma$ with $\mathscr{F}\cap \Gamma_{1}=\emptyset$, and if $\Gamma_{2}< \Gamma_{1}$ is a subgroup of finite index, not necessarily normal, then the proof of Hummel and Schroeder also shows that $X_{\Gamma_{2}}$ carries a metric of nonpositive curvature with the same properties as above. Consequently, if $b_{1}(\Gamma)$ is positive, the result of our main theorem will also apply to the filling $\overline{\Gamma_{2}}$. 
\end{rem}

We will say that a lattice $\Gamma_{1} < \PUD$ has {\it very nice cusps} if it has nice cusps as in Definition~\ref{nicecusps} and if $X_{\Gamma_{1}}$ carries a nonpositively curved Riemannian metric with all the properties from the paragraph above. Using the residual finiteness of lattices, it is easy to see that any nonuniform lattice in $\PUD$ has finite index subgroups which are lattices with very nice cusps. 

Let us say a word about the proof of these results: any horoball of $\HC$ is diffeomorphic to $(a,b]\times N$. Hummel and Schroeder consider in~\cite{hs} Riemannian metrics on $(a,b]\times N$ of the form $\langle \cdot , \cdot \rangle_{f,g}:=dt^{2}\oplus \mu_{f(t),g(t)}$ where $f$ and $g$ are smooth positive functions. Recall here that $\mu_{f(t),g(t)}$ is the left-invariant metric on $N$ obtained by rescaling $\mu$ by $f(t)^{2}$ on $\frak{z}$ and by $g(t)^{2}$ on $\frak{r}$. If $f$ and $g$ coincide with one of the model functions $\alpha e^{-t}$ and $\alpha^{2} e^{-2t}$ neat $t=b$, one can glue isometrically $(a,b]\times N$ endowed with this metric to the exterior of a horoball. This metric is also invariant by left translations on the $N$ factor, one can thus consider it as a metric on any (truncated) cusp of the form $(a,b]\times \Lambda \backslash N$ for a lattice $\Lambda < N$. One can then impose conditions on $f$ and $g$ first to guarantee that this metric extends smoothly to the compactified cusp and then to guarantee that it is either K\"ahler or nonpositively curved (but one cannot do both at the same time). We refer the reader to~\cite[\S 3]{hs} for more details.


The only consequences of Hummel and Schroeder's result that we will need is that if $\Gamma$ has very nice cusps, then $X_{\Gamma}$ is aspherical and the Abelian subgroups of its fundamental group are understood, as shown by the following proposition: 
 
\begin{prop}\label{flats} Let $\Gamma < \PUD$ be a lattice with very nice cusps. Let $i : \Z^{2}\to \pi_{1}(X_{\Gamma})\simeq \overline{\Gamma}$ be an injective homomorphism. Then $i(\Z^{2})$ is conjugated to a subgroup of the fundamental group of one of the totally geodesic elliptic curves in $X_{\Gamma}-\HC/\Gamma$.  
\end{prop}
{\it Proof.} Consider the action of $\overline{\Gamma}$ on the universal cover $\widetilde{X_{\Gamma}}$ of $X_{\Gamma}$, endowed with the lift of a nonpositively curved Riemannian metric as above. By~\cite[II.6.10]{bh}, every element of $\overline{\Gamma}$ is a semisimple isometry of $\widetilde{X_{\Gamma}}$. By the flat torus Theorem~\cite[II.7.1]{bh}, the group $i(\Z^{2})$ must preserve a totally geodesic flat $\R^{2} \hookrightarrow \widetilde{X_{\Gamma}}$. But since Hummel and Schroeder's Riemannian metric has negative curvature in the open set $\HC/\Gamma \subset X_{\Gamma}$, this flat must be contained in the inverse image of one of the elliptic curves from $X_{\Gamma}-\HC/\Gamma$. This gives the desired result.\hfill $\Box$

\begin{rem} The fundamental groups of complex surfaces such as $X_{\Gamma}$ fit into the study of groups acting on ${\rm CAT}(0)$ spaces with isolated flats or of relatively hyperbolic groups as in~\cite{groves,hk} or~\cite{osin} for instance. 
\end{rem}

The following result is well-known. 

\begin{prop} Let $\Gamma < \PUD$ be a nonuniform lattices with nice cusps. The K\"ahler surface $X_{\Gamma}$ constructed above is algebraic. 
\end{prop}
{\it Proof.} The reader will find a proof of this fact in~\cite[\S 2]{dicerbo}, based on the classification of surfaces. A different proof, applying also in higher dimensions, goes as follows. Mok~\cite{mok} proved that there exists a holomorphic map $X_{\Gamma} \to \mathbb{P}^{N}$ (for some $N$) which is an embedding on the open set $\HC/\Gamma \subset X_{\Gamma}$ and which contracts the elliptic curves in $X_{\Gamma}-\HC/\Gamma$ to points. This implies that $X_{\Gamma}$ is Moishezon. But a K\"ahler manifold which is also Moishezon is projective.\hfill $\Box$


\section{Coherence and homomorphisms to Abelian groups for K\"ahler groups}\label{rappelstrategie}

In this section we recall a result essentially due to Kapovich~\cite{kgt13}, which was used in his proof of the noncoherence of cocompact arithmetic lattices of the simplest type in ${\rm PU}(n,1)$. The result concerns fundamental groups of compact aspherical K\"ahler surfaces $X$ with positive first Betti number. Recall that for a K\"ahler manifold, the first Betti number is even. If $b_{1}(X)>0$, one can thus consider surjective homomorphisms from $\pi_{1}(X)$ to $\mathbb{Z}^{2}$. 

In the following we will say that a compact complex surface is a {\it Kodaira surface} if it admits a submersion onto a compact hyperbolic Riemann surface with connected hyperbolic fibers. Although this definition is not completely standard (see the discussion in~\cite{kotschick}), we will use it here. Such a surface is necessarily aspherical.

\begin{theorem}\label{kapovich-alt} (Kapovich) Let $X$ be an aspherical K\"ahler surface with positive first Betti number. Assume that $\pi_{1}(X)$ has no finitely generated Abelian subgroup whose normalizer has finite index in $\pi_{1}(X)$. Then, one of the following three cases occurs:
\begin{enumerate}
\item The group $\pi_{1}(X)$ is not coherent.
\item The surface $X$ has a finite cover which is a Kodaira surface.
\item For every surjective homomorphism $\phi : \pi_{1}(X) \to \mathbb{Z}^{2}$, the kernel of $\phi$ is isomorphic to the fundamental group of a closed Riemann surface. 
\end{enumerate}
\end{theorem}

The proof of this theorem has several ingredients that we now list. 

\begin{enumerate}
\item One of them is {\it Delzant's alternative}~\cite{d}, stating the following: if $X$ is a closed K\"ahler manifold and if $\phi : \pi_{1}(X)\to A$ is a homomorphism to an Abelian group, then the kernel of $\phi$ is finitely generated unless $X$ admits a holomorphic fibration onto a hyperbolic $2$-dimensional orbifold. We refer the reader to~\cite{d} for a more precise statement and for the definition of hyperbolic $2$-dimensional orbifolds; here we will only need this weak form of Delzant's Theorem.  
\item We will use the fact that Poincar\'e duality groups of dimension $2$ are fundamental groups of closed Riemann surfaces, as follows from the work of Eckmann together with Bieri, Linnel and M\"uller, see~\cite{eckmann} and the references there. For the definition of Poincar\'e duality groups, we refer the reader to~\cite[VIII.10]{browncoho}. 
\item If $G$ is a Poincar\'e duality group of dimension $4$ and if one has a a short exact sequence
\begin{center}$\begin{CD}
1 @>>> H @ >>> G @>>> \pi_{1}(S) @>>> 1 \\
\end{CD}$\end{center}
where $H$ is finitely presented and $S$ is a closed Riemann surface, then $H$ is a Poincar\'e duality group of dimension $2$. This result is due to Hillman, see Theorem 1.19 in~\cite{hillman2002} (which is more general). Combined with the previous result, Hillman's result implies that $H$ is the fundamental group of a closed Riemann surface.  
\item Let $X$ be an aspherical K\"ahler surface whose fundamental group fits into a short exact sequence
\begin{center}$\begin{CD}
1 @>>>  H_{1} @>>> \pi_{1}(X) @>\pi>> H_{2} @>>> 1 \\
\end{CD}$\end{center}
where both $H_{1}$ and $H_{2}$ are fundamental groups of closed Riemann surfaces of genus greater than one. Then the homomorphism $\pi$ is induced by a holomorphic submersion with connected fibers onto a closed Riemann surface. This result, or slight variations on it, has been proved independently by several people at the end of the 90's, see for instance~\cite{hillman2000,k1998,kotschick}. We refer the reader to~\cite{kotschick} for a short elegant proof.  
\end{enumerate}

We now turn to the proof of Theorem~\ref{kapovich-alt}, based on the above ingredients. Although this proof is essentially contained in~\cite{kgt13}, we will explain it for the reader's convenience. We start with the:

\begin{prop}\label{step} Let $X$ be an aspherical K\"ahler surface with coherent fundamental group. Let $\phi : \pi_{1}(X) \to \Z^{2}$ be a homomorphism with finitely generated kernel. Then, the kernel of $\phi$ is isomorphic to the fundamental group of a closed Riemann surface. 
\end{prop}
{\it Proof.} Since $\pi_{1}(X)$ is coherent the kernel of $\phi$ is finitely presented. Now by Hillman's theorem mentioned above, the kernel of $\phi$ must be a Poincar\'e duality group of dimension $2$ (for short: a ${\rm PD}(2)$ group). The characterization of ${\rm PD}(2)$ groups then implies that the kernel of $\phi$ is the fundamental group of a closed Riemann surface.\hfill $\Box$

{\it Proof of Theorem~\ref{kapovich-alt}.} We assume that the fundamental group of $X$ is coherent and prove that it must satisfies the second or the third possibility from the theorem. If every homomorphism $\pi_{1}(X)\to \Z^{2}$ has finitely generated kernel, then Proposition~\ref{step} implies that $\pi_{1}(X)$ satisfies the third possibility of the theorem. 

Now if there exists one homomorphism $\phi_{0} : \pi_{1}(X)\to \Z^{2}$ whose kernel is not finitely generated, then Delzant's Theorem implies that there is a holomorphic fibration $\pi : X \to \Sigma$ onto a $2$-dimensional hyperbolic orbifold. Such a fibration induces a surjective homomorphism
$$\pi_{1}(X)\to \pi_{1}^{orb}(\Sigma)$$
with finitely generated kernel, where $\pi_{1}^{orb}(\Sigma)$ is the orbifold fundamental group of $\Sigma$; note that we implicitly assume here that the orbifold structure on $\Sigma$ is given by the multiplicities of the singular fibers of $\pi$. See~\cite[\S 4.1]{delzant2008} for all of this. There exists a finite cover $X_{1}\to X$ and a finite orbifold cover $\Sigma_{1}\to \Sigma$ such that $\Sigma_{1}$ is a manifold and such that $\pi$ lifts to a holomorphic map $\pi_{1} : X_{1}\to \Sigma_{1}$ inducing a surjective homomorphism $$(\pi_{1})_{\ast} : \pi_{1}(X_{1})\to \pi_{1}(\Sigma_{1}).$$ 
The kernel of $(\pi_{1})_{\ast}$ is also finitely generated. Since $\pi_{1}(X)$ is assumed to be coherent (and since $\pi_{1}(X_{1})$ is a subgroup of $\pi_{1}(X)$), the kernel of $(\pi_{1})_{\ast}$ must be finitely presented. Using Hillman's result again, we obtain that $Ker((\pi_{1})_{\ast})$ is the fundamental group of a closed oriented surface $F$. The hypothesis on Abelian subgroups of $\pi_{1}(X)$ implies that the genus of $F$ is greater than 
$1$. Then the fourth result recalled above implies that $X_{1}$ is a Kodaira surface. Hence the surface $X$ satisfies the second possibility of the theorem. This concludes the proof.\hfill $\Box$


\section{Flats and Poincar\'e's theorem}\label{proof}

We now prove the main theorem. So let $\Gamma \subset {\rm PU}(2,1)$ be a torsion-free nonuniform lattice with nice cusps and such that $b_{1}(\Gamma)$ is positive. We have seen in section~\ref{rappelcusp} that there exists a finite set $\mathscr{F}$ of parabolic isometries of $\Gamma$ such that any finite index normal subgroup $\Gamma_{1}$ of $\Gamma$ with trivial intersection with $\mathscr{F}$ has very nice cusps, which means by definition that the compactified surface $X_{\Gamma_{1}}$ admits a nonpositively curved Riemannian metric enjoying all the properties described in section~\ref{rappelcusp}. We take exactly this set $\mathscr{F}$ as the set appearing in the statement of our theorem. We now fix a finite index normal subgroup $\Gamma_{1}$ of $\Gamma$ such that $\Gamma_{1}\cap \mathscr{F}=\emptyset$. We also endow once and for all the surface $X_{\sgr}$ with one of the nonpositively curved Riemannian metrics constructed by Hummel and Schroeder, whose properties were listed in Section~\ref{rappelcusp}. We start with the following proposition. 

\begin{prop}\label{autrestep} The group $\overline{\Gamma_{1}}$ does not contain any finitely generated Abelian group $A$ whose normalizer has finite index in $\overline{\Gamma_{1}}$. No finite cover of the surface $X_{\Gamma_{1}}$ is a Kodaira surface.  
\end{prop}
{\it Proof.} We denote by $\widetilde{X_{\sgr}}$ the universal cover of $X_{\sgr}$. Suppose that $A$ is a finitely generated Abelian group of $\overline{\sgr} \simeq \pi_{1}(X_{\sgr})$. For $g\in \overline{\Gamma_{1}}$, considered as an isometry of $\widetilde{X_{\Gamma_{1}}}$, let $${\rm Min}(g)=\{x\in \widetilde{X_{\Gamma_{1}}}, d(x,g(x))={\rm min}(g)\}$$ where ${\rm min}(g)$ is the translation length of $g$. Let ${\rm Min}(A)=\cap_{g\in A}{\rm Min}(g)$. It is known that ${\rm Min}(A)$ is a convex subset of $\widetilde{X_{\Gamma_{1}}}$ which splits as a product $Y\times \R^{r}$ where $r$ is the rank of $A$, see~\cite[II.7]{bh}. If $r\ge 2$ or if $Y$ is not reduced to a point, one sees, using that the curvature is negative on the open set $\HC/\Gamma_{1}\subset X_{\Gamma_{1}}$, that ${\rm Min}(A)$ must be contained in a connected component $\widetilde{E}$ of the inverse image of an elliptic curve $E\subset X_{\Gamma_{1}}$. Since the normalizer $N(A)$ of $A$ in $\overline{\sgr}$ preserves ${\rm Min}(A)$ by~\cite[II.7]{bh}, this actually implies that $N(A)$ must preserve the flat $\widetilde{E}$. Any element $g$ of $N(A)$ acts on $\widetilde{E}$ as a translation (being semisimple and orientation preserving). This implies that $N(A)$ is free Abelian of rank at most $2$, hence it cannot be of finite index in $\overline{\sgr}$. If $r=1$ and $Y$ is a point, the proof is similar (and simpler): ${\rm Min}(A)$ is then made of a single geodesic, which must be $N(A)$-invariant, preventing $N(A)$ from being of finite index. 

For the second claim of the proposition, suppose that a finite cover $X_{1}$ of $X_{\sgr}$ is a Kodaira surface, i.e. that there exists a closed hyperbolic Riemann surface $S$ and a holomorphic submersion $\pi : X_{1} \to S$ with connected hyperbolic fibers. Let $E\subset X_{1}$ be a totally geodesic elliptic curve. Note that we have seen that such curves exist in $X_{\sgr}$; hence they also exist in any finite cover of $X_{\sgr}$. Since $S$ is Kobayashi hyperbolic, the restriction of $\pi$ to $E$ is constant; hence $E\subset \pi^{-1}(\ast)$ for some point $\ast \in S$. Since $\pi^{-1}({\ast})$ is smooth and connected, this implies that $E=\pi^{-1}(\ast)$. This contradicts the fact that $\pi$ has hyperbolic fibers.\hfill $\Box$

To prove that $\overline{\Gamma_{1}}$ is not coherent, we will now apply Theorem~\ref{kapovich-alt}. Note that $\overline{\Gamma_{1}}$ indeed satisfies the hypothesis of that theorem: the first Betti number $b_{1}(\Gamma)$ of $\Gamma$ was assumed to be positive hence $b_{1}(\Gamma_{1})$ is also positive; moreover $\Gamma_{1}$ and $\overline{\Gamma_{1}}$ have the same first Betti number, hence $b_{1}(\overline{\Gamma_{1}})>0$. Also, by Proposition~\ref{autrestep}, the group $\overline{\Gamma_{1}}$ satisfies the hypothesis on Abelian subgroups appearing in Theorem~\ref{kapovich-alt}. 

Observe that $\overline{\Gamma_{1}}$ cannot satisfy the second possibility appearing in Theorem~\ref{kapovich-alt}, thanks to Proposition~\ref{autrestep}. So we must exclude the third possibility of that theorem, this will imply that $\overline{\Gamma_{1}}$ is not coherent. We now assume by contradiction that we are in the third case given by Theorem~\ref{kapovich-alt}. Pick any surjective homomorphism $f : \ssg \to \Z^{2}$. We know that there is a short exact sequence 
\begin{center}$\begin{CD}
1 @>>> \pi_{1}(S) @>>> \ssg @>f>> \mathbb{Z}^{2} @>>> 1 \\
\end{CD}$\end{center}
where $S$ is a closed Riemann surface. The surface $S$ is necessarily of genus greater than $1$ by Proposition~\ref{autrestep}. This short exact sequence determines a homomorphism $$\Psi : \mathbb{Z}^{2} \to Out(\pi_{1}(S)).$$
Here $Out(\pi_{1}(S))$ is the group of outer automorphisms of $\pi_{1}(S)$, also known as the mapping class group of $S$. As in~\cite{kgt13}, we will use the fact that a rank $2$ Abelian subgroup of the mapping class group of $S$ must contain a reducible element, see~\cite{blm}. This implies that there exists an embedding $i : \mathbb{Z}^{2} \to \ssg$ such that the composition of $i$ with the projection $f : \ssg \to \mathbb{Z}^{2}$ has rank one image. Indeed if $\Psi$ is not faithful, the existence of such an embedding $i$ is clear; if $\Psi$ is faithful, this follows from the above mentioned fact. But according to Proposition~\ref{flats}, and up to conjugacy, any embedding $i : \mathbb{Z}^{2}\to \ssg$ has its image contained into the image of the map 
$$\pi_{1}(E)\to \pi_{1}(X_{\Gamma_{1}})=\ssg$$
induced by the inclusion of a totally geodesic elliptic curve $E$ in $X_{\Gamma_{1}}$. This implies:

\begin{prop}\label{sup} If the third possibility of Theorem~\ref{kapovich-alt} occurs, there exists a totally geodesic elliptic curve $E\hookrightarrow X_{\Gamma_{1}}$ such that the holomorphic map 
$$h : E \to {\rm Alb}(X_{\Gamma_{1}})$$
obtained by composition of the inclusion of $E$ in $X_{\Gamma_{1}}$ and the Albanese map of $X_{\Gamma_{1}}$ is nontrivial. 
\end{prop}

Recall that the Albanese variety ${\rm Alb}(Y)$ of a compact K\"ahler manifold $Y$ is a compact complex torus of complex dimension equal to $b_{1}(Y)$ endowed with a holomorphic map ${\rm alb}(Y) : Y \to {\rm Alb}(Y)$ which induces an isomorphism
$${\rm alb}(Y)^{\ast} : H^{1}({\rm Alb}(Y),\R)\to H^{1}(Y,\R)$$
between the first cohomology groups of $Y$ and ${\rm Alb}(Y)$. The map ${\rm alb}(Y)$ is canonical up to translation. The torus ${\rm Alb}(Y)$ is algebraic if $Y$ is (see~\cite{voisin}). In the proposition above, we can assume that $h(E)$ is a subtorus of ${\rm Alb}(X_{\Gamma_{1}})$, up to composing the map ${\rm alb}(X_{\Gamma_{1}})$ with a translation. We will make this assumption below. 

We now continue the proof. Since the surface $X_{\Gamma_{1}}$ is algebraic, so is ${\rm Alb}(X_{\Gamma_{1}})$. We will apply Poincar\'e's reducibility theorem to ${\rm Alb}(X_{\Gamma_{1}})$. Recall that this theorem states that if $A$ is an Abelian variety and $B$ a subtorus of $A$, there exists another subtorus $C\subset A$ such that there is an isogeny $B\times C \to A$; see~\cite[VI.8]{debarre} for a proof. We apply Poincar\'e's theorem to $A={\rm Alb}(X_{\Gamma_{1}})$ and $B=h(E)$ the image of the elliptic curve appearing in Proposition~\ref{sup}. Hence there exists an Abelian variety $C$ of dimension one less than ${\rm Alb}(X_{\Gamma_{1}})$ and an isogeny $h(E)\times C \to {\rm Alb}(X_{\Gamma_{1}})$ which itself gives rise to an isogeny $$u : E \times C \to {\rm Alb}(X_{\Gamma_{1}}).$$ 
Let $\pi : X_{2}\to X_{\Gamma_{1}}$ be a finite cover of $X_{\Gamma_{1}}$ such that the map ${\rm alb}(X_{\sgr}) \circ \pi : X_{2} \to {\rm Alb}(X_{\sgr})$ lifts to $E \times B$ i.e. such that there exists a holomorphic map $g : X_{2}\to E \times B$ making the following diagram commutative:

\begin{center}$\begin{CD}
X_{2} @>g>> E \times B \\
@VV\pi V @VVuV \\
X_{\sgr} @>alb(X_{\sgr})>> {\rm Alb}(X_{\sgr}) \\
\end{CD}$\end{center}
Denote by $g_{2}$ the composition of $g$ with the first projection from $E \times B$ to $E$. We are now going to repeat the previous line of arguments, but considering the homomorphism $$(g_{2})_{\ast} : \pi_{1}(X_{2}) \to \pi_{1}(E)\simeq \Z^{2}$$
induced by $g_{2}$. The main difference is that we now know that this homomorphism is induced by a holomorphic map. We apply again Theorem~\ref{kapovich-alt}, this time to $X_{2}$ instead of $X_{\sgr}$. If  $\pi_{1}(X_{2})$ is not coherent, the same is true for $\overline{\sgr}$. So we must exclude the second and third possibilities of the theorem for $X_{2}$. The fact that $X_{2}$ has not finite cover which is a Kodaira surface follows from Proposition~\ref{autrestep}. Assume now that the third possibility of Theorem~\ref{kapovich-alt} holds and consider the homomorphism $(g_{2})_{\ast}$. Its kernel is then a non-Abelian surface group. As before the image of the homomorphism $\Z^{2}\to Out(Ker((g_{2})_{\ast})$ induced by $g_{2}$ must have rank $<2$ or must contain a reducible element, thus giving rise to a rank $2$ Abelian subgroup $$A\simeq \Z^{2} < \pi_{1}(X_{2})$$ such that the restriction of $(g_{2})_{\ast}$ to $A$ has rank $1$ image and rank $1$ kernel. But now we can assume again that, up to conjugacy, $A$ is contained in the fundamental group of a certain elliptic curve $i_{1} : E_{1}\hookrightarrow X_{2}$. We thus obtain a homomorphism $$(g_{2}\circ i_{1})_{\ast} : \pi_{1}(E_{1})\to \pi_{1}(E)$$
between the fundamental groups of two elliptic curves which has rank $1$ image and is induced by a holomorphic map. This is a contradiction. Hence the groups $\pi_{1}(X_{2})$ and $\overline{\Gamma_{1}}$ are not coherent. This completes the proof of our main theorem.

Let us make one final remark, valid in any dimension,  concerning nonuniform lattices in ${\rm PU}(n,1)$, and their fillings. If $\Gamma < {\rm PU}(n,1)$ is a torsionfree nonuniform lattice with nice cusps (i.e. whose parabolic elements have no rotational part) the positivity of the first Betti number of $\Gamma$ is equivalent to that of $b_{1}(X_{\Gamma})$. Similarly, the group $\Gamma$ surjects onto a non-Abelian free groups if and only if $\pi_{1}(X_{\Gamma})$ does. This should motivate the study of the spaces $X_{\Gamma}$ and of their fundamental groups, since the positivity of the first Betti number as well as the largeness of lattices in ${\rm PU}(n,1)$ are well-known open problems. 


\bigskip
\bigskip

\begin{small}
\begin{tabular}{l}
IRMA, Universit\'e de Strasbourg \& CNRS\\
67084 Strasbourg, France\\
ppy@math.unistra.fr\\    
\end{tabular}
\end{small}

 \end{document}